\documentclass[12pt]{article}
\usepackage{amsmath,amsthm,amsfonts,amscd,amssymb,latexsym,eucal}

\def\qed{{\unskip\nobreak\hfil\penalty50
\hskip2em\hbox{}\nobreak\hfil$\square$
\parfillskip=0pt \finalhyphendemerits=0\par}\medskip}
\def\proof{\trivlist \item[\hskip \labelsep{\bf Proof.}]}
\def\endproof{\null\hfill\qed\endtrivlist}

\def\a{\alpha}

\def\de{\delta}

\def\ep{\epsilon}
\def\f{\varphi}

\def\i{\iota}

\def\phi{\varphi}

\def\oma{\otimes_{\max}}
\def\omi{\otimes_{\min}}

\newtheorem{theorem}{Theorem}
\newtheorem{lemma}[theorem]{Lemma}

\newtheorem{corollary}[theorem]{Corollary}

\newtheorem{proposition}[theorem]{Proposition}

\def\2#1{{\cal #1}}

\def\A{{\cal A}}

\def\E{{\cal E}}
\def\H{{\cal H}}
\def\K{{\cal K}}

\def\C{{\mathbb C}}

\def\t#1{\text{\rm #1}}
\def\d{\text{\rm d}}

\title{\bf A Remark on\\
Quantum Group Actions\\  
and Nuclearity}

\author{
{\sc S. Doplicher}\\
Dipartimento di Matematica\\
Universit\`a di Roma ``La Sapienza''\\
Piazzale A. Moro 5, I-00185 Roma, Italy\\
\vphantom{X}\\
{\sc R. Longo, J.E. Roberts, L. Zsid\'o}\\
Dipartimento di Matematica\\
Universit\`a di Roma ``Tor Vergata''\\
Via della Ricerca Scientifica 1, I-00133 Roma, Italy}
\date{}
\begin{document}

\maketitle

\centerline{{\large {\sl Dedicated to Huzihiro Araki}}} 
 \centerline{{\large {\sl on the occasion of his seventieth birthday}}}
\bigskip

\begin{abstract} Let $H$ be a compact quantum group with faithful Haar 
measure and bounded counit. If $\a$ is an action of $H$ on a 
$C^*$-algebra $A\,$, we show that $A$ is nuclear if and only if the fixed-point 
subalgebra $A^{\a}$ is nuclear. As a consequence $H$ is a nuclear $C^*$-algebra.
\end{abstract}
\vfill

\thanks{{\footnotesize Supported in part by MIUR and GNAMPA-INDAM.}}

\thanks{{\footnotesize E-mails: dopliche@mat.uniroma1.it, longo@mat.uniroma2.it,  
roberts@mat.uniroma2.it, zsido@mat.uniroma2.it}}

\newpage

\section{Introduction}
The notion of compact quantum group has been axiomatized by 
Woronowicz \cite{W2} and some variations appeared in subsequent papers, 
e.g. \cite{BS,V,VW,KW}. Here we shall consider compact quantum groups
allowing a faithful Haar state and bounded counit, conditions that cover
a rather wide range of applicability cf. \cite{W1}.

In this note we will show that a $C^*$-algebra acted upon by a compact 
quantum group with faithful Haar state and bounded counit is nuclear
if and only if the fixed point algebra is nuclear. In particular, a
compact quantum group with faithful Haar state and bounded counit is
a nuclear $C^*$-algebra. Notice however that our assumptions about
the faithfulness of the Haar state and the boundedness of the counit
are essential: there exist non-nuclear compact quantum groups, as we
can infer from \cite{W1}, second remark after the statement of
Proposition 1.8, and \cite{VW}.

As a special case, a $C^*$-algebra acted upon by a compact group
of automorphisms is nuclear if and only if the fixed point algebra is
nuclear. Indeed this was our initial motivation, during a common
discussion long ago at the Rome seminar on Operator Algebras, inspired 
by a result of H\o egh-Krohn, Landstad and St\o rmer \cite{HLS} showing
that there is a compact ergodic action on a $C^*$-algebra only if the
algebra is nuclear, and by a lemma in \cite{DR}. 
This particular case is discussed separately as it is a simple and 
clarifying illustration for the general case. 

We notice that another extension of the result of H\o egh-Krohn,
Landstad and St\o rmer to compact matrix pseudogroups was already
proved by Boca in \cite{B}, Corollary 23 : any $C^*$-algebra acted
upon ergodically by a nuclear compact matrix pseudogroup is again
nuclear. Actually Boca's proof works in the case of an ergodic action
of an arbitrary nuclear compact quantum group.

We also analyse compact actions on von Neumann algebras where the
analogous result holds for injective von Neumann algebras. This 
$W^*$-version can however also be inferred from a similar result for 
crossed products due to Connes \cite{C2} and indeed can be extended
in part to the case of integrable actions of locally compact groups
or Kac algebras, but we do not treat this case here.

Our discussion is elementary, but in the last section ($W^*$-case) 
we make use of the relation between injectivity and semi-discreteness 
due to Effros-Lance, Choi-Effros and Connes \cite{EL,CE,C1}.
\medskip

\noindent
{\bf Notation:}
\smallskip

If $\A$ is a ${}^*$-algebra we shall denote  the identity map on $\A$
by $\i_{\substack{{}\\ \A}}$, or simply by $\i$. When $\A$
is unital, we shall denote its unit by $1_\A\,$.

The maximal tensor product for $C^*$-algebras is denoted by $\oma\,$,
the minimal one by $\omi$ or simply by $\otimes\,$, the algebraic
tensor product by $\odot\,$. For a $C^*$-algebra $A$, we shall usually
identify $\C\otimes A$ and $A\otimes \C$ with $A\,$.

By a homomorphism we shall always mean a $^*$-homomorphism.

\section{Basic Result}
Let us call a conditional expectation $\E$ on a $C^*$-algebra $A$
{\it GNS-faithful} if
$$
x\in A\;\; \&\;\; \E (y^*x^*xy)=0\;\text{ for }\;
\forall y\in A\;\implies\; x=0\ .
$$
This means that the direct sum of the GNS representations associated 
with all $\E$-invariant states is injective. Of course, $\E$ is 
GNS-faithful if $\E$ is faithful.

Let $A$ be a $C^*$-algebra, $A_0\subset A$ a $C^*$-subalgebra and 
$\E: A\to A_0$ a conditional expectation.

If $B$ is any $C^*$-algebra, it is easy to see that
$\i_B\odot\E:B\odot A\to B\odot A_0$ extends to a conditional
expectation from $B\omi A$ to $B\omi A_0\,$, which is
(GNS-)faithful if $\E$ is (GNS-)faithful and which we
denote by $\i_B\otimes\E\,$.

We shall say that $\E$ is {\it stably (GNS-)faithful} if, for
every $C^*$-algebra $B$, $\iota_B\odot\E$ extends to a bounded
(GNS-)faithful map  $\tilde\E: B\oma A\to B\oma A$. In this
case $\tilde\E$ is a (GNS-)faithful conditional expectation from
$B\oma A$ onto the closure of $B\odot A_0$ in $B\oma A$.
\begin{proposition}\label{stably}
Let $A$ be a $C^*$-algebra, $A_0\subset A$ a $C^*$-subalgebra and 
$\E: A\to A_0$ a stably GNS-faithful conditional expectation.
Then $A$ is nuclear iff $A_0$ is nuclear.
\end{proposition}
\proof 
The implication $A$ nuclear $\implies$ $A_0$ nuclear is known (indeed 
one just needs $A_0$ to be the range of a conditional expectation). 
One way to see this is to recall that $A$ is nuclear iff the enveloping 
von Neumann algebra $A^{**}$ is injective \cite{EL,CE}, and to consider 
the double transposed conditional
expectation $\E^{**}:A^{**}\to A_0^{**}$. If $A$ is nuclear, $A^{**}$ is
injective, and then $A_0^{**}$ is injective too, hence $A_0$ is nuclear.

Conversely, assume that $A_0$ is nuclear. Given a $C^*$-algebra $B$,
the identity map on the algebraic tensor product $B\odot A$ extends
to a homomorphism
$$
\pi:B\otimes_{\max} A\to B\otimes_{\min} A\ .
$$
To show that $A$ is nuclear, we have to prove $\pi$ to be one-to-one. 

By assumption, there is a GNS-faithful conditional 
expectation $\tilde\E:B\otimes_{\max}A\to B\otimes_{\max}A$ extending
$\iota_B\odot\E$. 
Denoting the conditional expectation $\i_B\otimes\E$ on $B\omi A$ by
$\tilde\E'$, we clearly have a commutative diagram
\begin{equation}\label{cd}
\left.
\begin{CD}
B\otimes_{\max} A   @> \pi >> B\otimes_{\min} A 
\\ @V\tilde\E V V       @VV\tilde\E' V  \\  B\otimes_{\max} A  
@>\pi >>  B\otimes_{\min} A \end{CD} \right.\quad .
\end{equation}
Since $\tilde\E$ maps $B\odot A$ onto $B\odot A_0$, by continuity 
$\tilde\E$ maps $B\otimes_{\max} A$ onto the closure of $B\odot A_0$ in 
$B\otimes_{\max} A$, that we may denote by $B\otimes A_0$ since $A_0$
is nuclear. Then (\ref{cd}) yields the commutative diagram  
\begin{equation}\label{cd2}
\left.
\begin{CD}
B\otimes_{\max} A   @> \pi >> B\otimes_{\min} A 
\\ @V\tilde\E V V       @VV\tilde\E' V  \\  B\otimes A_0  
@>\i_{\substack{{}\\ B\otimes A_0}} >>  B\otimes A_0 \end{CD}
\right.\quad .
\end{equation}

Now if $x\in B\oma A$ belongs to the kernel of $\pi$, then 
$\tilde\E'\big(\pi(y^*x^*xy)\big)=0$ for all $y\in B\oma A\,$. By the
commutativity of diagram (\ref{cd2}) we have $\tilde\E(y^*x^*xy)=0$
for all $y\in B\oma A\,$, so $x=0$ because $\tilde\E$ is GNS-faithful.
We conclude that $\pi$ is injective and $A$ nuclear.
\endproof

\section{Compact Group Actions and Nuclearity}
In the sequel, group actions on $C^*$-algebras and
von Neumann algebras will be assumed continuous in the usual sense, 
namely pointwise norm continuity in the $C^*$-case and
pointwise weak$^*$-continuity in the $W^*$-case.

\begin{proposition}\label{groupact} 
Let $\alpha:G\to\t{Aut}(A)$ be an action
of a compact group $G$ on a $C^*$-algebra $A$. Then $A$ is nuclear
iff the fixed point $C^*$-subalgebra $A^\alpha$ is nuclear.
\end{proposition}

\proof Let $\E_\alpha:A\to A^\alpha$ be the conditional expectation 
defined by
$$
\E_\alpha(a)=\int\alpha_g(a)\d g\ ,\quad a\in A\ .
$$
According to Proposition \ref{stably}, it suffices to show that
$\E_\alpha$ is stably faithful.

The action $\beta=\i_B\otimes\alpha$ on $B\odot A$ preserves 
the maximal cross norm $\Vert\cdot\Vert_{\max}$. Furthermore, $G\ni g
\mapsto \beta_g(x)$ is continuous for every $x\in B\odot A$ with respect
to $\Vert\cdot\Vert_{\max}$.
Thus $\beta$ extends to an action $\beta^{\max}$ of $G$ on $B\oma A$.

Set 
$$
\E_\beta(x)=\int\beta^{\max}_g(x)\,\d g\, ,\quad x\in B\oma A\ .
$$
Then $\E_{\beta}$ is a faithful conditional expectation because  we have
\[
\E_\beta(x)=\int\beta^{\max}_g(x)\,\d g=0\;\implies\;\beta^{\max}_g(x)=0 
\: \forall g\in G
\] 
for  every positive $x\in B\oma A\,$. 
Since $\E_\beta$ maps $B\odot A$ onto $B\odot A^\alpha$ and 
$\E_\beta\, |\, _{B\odot A}= \i\odot\E_{\alpha}$, we see that 
$\E_{\alpha}$ is stably faithful.
\endproof

\section{Quantum Group Actions and Nuclearity}
Following Woronowicz and Van Daele \cite{W2,V}, by a
{\it compact quantum group} we shall mean a unital $C^*$-algebra $H$ equipped with
a {\it comultiplication} $\de$, that is a unital homomorphism
$\de:H\to H\otimes H$ satisfying
\[
(\de\otimes\iota_H )\circ\de=(\iota_H\otimes\de )\circ\de\, ,
\]
such that
\[
\de(H)(H\otimes 1_H)\;\text{ and }\;\de(H)(1_H\otimes H)\;
\text{ are linearly dense in }\; H\otimes H\, .
\]
According to \cite{W2,V}, there then exists a unique {\it Haar state}
on $H\,$, that is a state $\f$ which satisfies the invariance condition
\begin{equation}\label{inv}
\big((\f\otimes\iota_H )\circ\de\big)(x)=\big((\iota_H\otimes\f )
\circ\de\big)(x)=\f(x)\, 1_H\, ,\quad x\in H\, . 
\end{equation}
We notice that by condition (\ref{inv}) the fixed point algebra
\begin{equation}\label{erg}
H^{\de} = \{x\in H\, ;\, \de (x)=x\otimes 1_H\}\;\text{ is equal to }
\; \C\, 1_H\, .
\end{equation}

The unique Haar state $\f$ is not necessarily faithful. However,
the cyclic vector $\xi_\f$ of the GNS-representation $\pi_\f$ is
also separating for $\pi_\f (H)''\,$. Indeed, $\de$ lifts to a
comultiplication $\de_\f$ on $M_\f =\pi_\f (H)''$
(see e.g. \cite{Wa}, Theorem 2.4) and the vector state
$\omega_{\xi_\f}$ on $M_\f$ satisfies the invariance conditions
corresponding to (\ref{inv}), in particular ${M_\f}^{\de_\f} = \C\,
1_{M_\f}\,$. Now \cite{SVZ}, Lemma 0.2.4, implies that the support of
$\omega_{\xi_\f} | M_\f$ belongs to ${M_\f}^{\de_\f} = \C\, 1_{M_\f}\,$,
hence it is equal to $1_{M_\f}\,$.
Consequently, replacing $H$ with $H/\t{ker}(\pi_\f)$ if necessary,
one can consider the case where $\f$ is faithful, cf. e.g. \cite{W1},
remarks after Theorem 5.6, \cite{MN,Wa}.

According to \cite{W2}, Theorem 2.2, the linear span $\A$ of all
matrix elements of all finite-dimensional unitary representations of
$H$ is a dense ${}^*$-subalgebra of $H$ with $1_H\in\A\, ,\,\de (\A )
\subset \A \odot \A$ and there are unique linear maps
$\ep : \A\to\mathbb C$ and $\kappa : \A\to \A\,$, called respectively
{\it counit} and {\it coinverse} or {\it antipode}, such that
\begin{align*}
(\ep\odot\i_\A )\big(\de (a)\big) &= (\i_\A\odot\ep )\big(\de (a)\big)
= a\, ,\hspace{2.7 cm} a\in\A\, , \\
\big( m\circ (\kappa\odot\i_\A )\big)\big(\de (a)\big) &
= \big( m\circ (\i_\A \odot\kappa )\big)\big(\de (a)\big)
= \ep (a)\, 1_H\, ,\quad a\in\A\, ,
\end{align*}
where the linear map $m : \A\odot\A \to \A$ is defined by
$m (a\odot b) = ab\, ,\, a,b\in\A\,$.

The counit $\ep$ is a multiplicative positive linear functional
on $\A$ with $\ep (1_H)=1\,$, but it is in general not bounded.
However, it is often bounded and then it extends to a multiplicative
state on $H\,$, still denoted by $\ep\,$, which satisfies
\begin{equation}\label{counit}
(\ep\otimes\i_H )\circ\de = (\i_H\otimes\ep )\circ\de =\iota_H\, .
\end{equation}
We notice that (\ref{counit}) implies the injectivity of $\de\,$.

In the rest of this section $H$ will denote a compact quantum
group with faithful Haar state and bounded counit. This is the case
for compact groups and for the quantum $U(N)$-group, but not for
the dual of a non-amenable discrete group, see \cite{W1}, second
remark after the statement of Proposition 1.8. 

By a {\it coaction} $\a$ of $H$ on a $C^*$-algebra $A$ we mean 
a homomorphism $\a:A\to A\otimes H$ such that
\begin{align*}\label{coaction}
(\a\otimes\i_H )\circ\a = &\; (\i_A\otimes\de )\circ\a \\ 
(\i_A\otimes\ep )\circ\a = &\; \i_A\ .
\end{align*}
The last equation implies that $\a$ is injective.

The fixed-point subalgebra is then defined by
\[ A^\a = \{a\in A\, ;\, \a(a)=a\otimes 1_H\}\ . \]
Denoting 
$$
\E = (\iota_A\otimes\f )\circ\a : A\to A\ ,
$$ 
we have the known fact:
\begin{lemma}\label{faithful}
 $\E$ is a faithful conditional expectation from 
$A$ to $A^\a$.
\end{lemma}
\proof
$\E$ is a faithful map, being the composition of faithful 
maps. Clearly $A^\a$ is contained in the range of $\E$.

We now apply standard calculations and check that
$\E$ is idempotent. Indeed, identifying $A$ with
$A\otimes \C\,$, we get 
\begin{align*}
\E^2 
=&\,\big((\iota_A\otimes\f )\circ\a\big)
\circ(\iota_A\otimes\f)\circ\a \\
=&\, (\iota_A\otimes\f)\circ(\a\otimes\f)\circ\a \\
=&\, (\iota_A\otimes\f\otimes\f)\circ(\a\otimes\iota_H)\circ\a \\
=&\, (\iota_A\otimes\f\otimes\f)\circ(\iota_A\otimes\de)\circ\a\\
=&\,\Big(\iota_A\otimes\big(
\underbrace{(\f\otimes\f )\circ\de}_{=\,\f(\cdot)\cdot 1_\C\otimes 1_\C}
\big)\Big)\circ\a \\
=&\, (\iota_A\otimes\f)\circ\a=
\E\ .
\end{align*}
To see that $\E(A)\subset A^\a$ we compute:
\begin{align*}
\a\circ\E 
=&\, \a\circ(\iota_A\otimes\f)\circ\a \\
=&\, (\iota_A\otimes\i_H\otimes\f)\circ(\a\otimes\iota_H)\circ\a \\
=&\, (\iota_A\otimes\i_H\otimes\f)\circ(\iota_A\otimes\de)\circ\a \\
=&\,\Big(\iota_A\otimes\big(
\underbrace{(\i_H\otimes\f)\circ\de}_{=\, \f(\cdot)\cdot 1_H\otimes 1_\C}
\big)\Big)\circ\a \\
=&\, \big(\i_A\otimes (\f(\cdot)\cdot 1_H )\big)\circ\a
= \E\otimes 1_H \, .
\end{align*}
The rest is now clear.
\endproof
\begin{lemma}\label{rho}
Let $B,A$ and $H$ be $C^*$-algebras. The identity map on $B\odot 
A\odot H$ extends to a homomorphism
\[
\rho: B\oma (A\omi H)\to (B\oma A)\omi H.
\]
\end{lemma}
\proof
Let $B$ and $A$ act faithfully on a Hilbert space $\H$ in such a way 
that the $C^*$-algebra 
generated by $B$ and $A$ is $B\oma A$ and let $H$ act faithfully 
on a Hilbert space $\K$. The $C^*$-algebra 
generated by $B$, $A$ and $H$ on $\H\otimes\K$ is clearly $(B\oma A)\omi 
H$ and the $C^*$-algebra generated by $A$ and $H$ is $A\omi H$. Thus
$(B\oma A)\omi H$ contains commuting copies of $B$ and $A\omi H$. By 
the universal property of $\oma$, we have a natural homomorphism
$\rho: B\oma (A\omi H)\to (B\oma A)\omi H$. 
\endproof 
Let $\a:A\to A\otimes H$ be a coaction as above and $B$ a 
$C^*$-algebra. By the universality property of $\oma$, the map 
$\iota_B\odot \a: B\odot A\to B\odot (A\otimes H)$ extends to a 
homomorphism $\iota_B\otimes \a: B\oma A\to B\oma (A\otimes H)$. By 
composing it with the map $\rho$ in Lemma \ref{rho}, we get a 
homomorphism
\[
\tilde\a\equiv \rho\circ(\iota_B\otimes \a): B\oma A\to (B\oma A)\otimes 
H\ .
\]
\begin{lemma}
$\tilde\a$ is a coaction of $H$ on $B\oma A$.
\end{lemma}
\proof
We first check the ``counit'' condition 
$$
(\i_{\substack{{}\\ B\oma A}}\otimes\ep)\circ\tilde\a =
\i_{\substack{{}\\ B\oma A}} .
$$
Clearly we have 
\begin{equation*}
\big( (\i_{\substack{{}\\ B\oma A}}\otimes\ep)\circ\tilde\a\big)
(b\otimes 1_A) = (\i_{\substack{{}\\ B\oma A}}\otimes\ep)
(b\otimes 1_A\otimes 1_H) = b\otimes 1_A, \quad b\in B\ .
\end{equation*}
Note that $(\i_{\substack{{}\\ B\oma A}}\otimes\ep)\circ\rho$ is the
bounded map from $B\oma(A\otimes H)$ to $B\oma A\,$, $b\otimes a\otimes h
\mapsto b\otimes \big(\ep(h)\cdot a\big)\,$, $a\in A\,$, $b\in B\,$,
$h\in H\,$, thus
\[
(\i_{\substack{{}\\ B\oma A}}\otimes\ep)\circ\rho: b\otimes x
\mapsto b\otimes \big((\i_A\otimes\ep)(x)\big),\quad b\in B,\ x
\in A\otimes H\ .
\]
By the counit property of $\ep$, we then have
\begin{align*}
\big( (\i_{\substack{{}\\ B\oma A}}\otimes\ep)\circ\tilde\a\big)
(1_B\otimes a)=&\, \big((\i_{\substack{{}\\ B\oma A}}\otimes\ep)\circ
\rho\big) (1_B\otimes \a(a)) \\
=&\, 1_B\otimes\big((\i_A\otimes\ep) \a(a)\big) =1_B\otimes a\, ,
\quad a\in A.
\end{align*}
As $\ep$ is multiplicative $(\i_{\substack{{}\\ B\oma A}})\otimes\ep\circ
\tilde\a$ acts identically on $B\odot A\,$, hence on $B\oma A$ by continuity.
As a consequence, $\tilde\a$ is injective.

Since $\tilde\a$ is injective, the restriction of $\rho$ to
$(\i_B\otimes\a)(B\oma A)$ is 
injective too and we can define the following commuting diagrams
\begin{equation*}
\begin{CD}
B\oma A    @> \i_B\otimes\a >> B\oma(A\otimes H) @>\i_B\otimes(\a\otimes\i_H)
>\i_B\otimes(\i_A\otimes\de)>B\oma((A\otimes H)\otimes H)
\\ @VV \i_{B\oma A} V   @VV\rho V @VV\rho' V \\  B\oma A  
@>\tilde\a >>  (B\oma A)\otimes H @>\tilde\a\otimes\i_H>
\i_{B\oma A}\otimes\de>(B\oma A)\otimes H\otimes H
\end{CD}
\end{equation*}
where $\rho'$ is the natural map constructed by Lemma \ref{rho}, 
showing that $\tilde\a$ is a coaction.
\endproof
\begin{proposition}\label{esf} 
$\E$ is stably faithful.
\end{proposition}
\proof
Let $B$ a $C^*$-algebra. The conditional expectation $\tilde\E$ on 
$B\oma A$ extending $\i_B\otimes\E$ on $B\odot A$ is 
the one associated with the coaction $\tilde\a$, hence it is faithful by 
Lemma \ref{faithful}.
\endproof
\begin{corollary} Let $\a$ be a coaction of a compact quantum group  
$H$ with faithful Haar state and bounded counit on a $C^*$-algebra $A$.
Then $A$ is nuclear iff $A^\a$ is nuclear.
\end{corollary}
\proof
This is a consequence of Propositions \ref{stably} and \ref{esf}.
\endproof
\begin{corollary} A compact quantum group $H$ with faithful Haar state
and bounded counit is a nuclear $C^*$-algebra.
\end{corollary}
\proof
The comultiplication $\de$ is a coaction of $H$ on itself and the
statement follows by the above corollary, taking (\ref{erg}) into account.
\endproof

\section{Compact Group Actions: the $W^*$-case}
Before giving a $W^*$-version of Proposition \ref{groupact}, 
we need some preliminaries.

Let $X$ be a Banach space, and $Y$ a weak$^*$ dense Banach subspace of
$X^*\,$. By a $\sigma (X,Y)$-continuous group of isometries $V:G\to B(X)$
we mean a $\sigma(X,Y)$-continuous homomorphism of the group $G$ into
the group of all $\sigma(X,Y)$-continuous linear isometries of $X$. 
We note that to check the $\sigma (X,Y)$-continuity of the map
$g\in G\to V_g\in B(X)$ we may verify that the maps $g \in G\to y(V_gx)$
are continuous when $x$ varies in a norm dense subset of $X$ and $y$
varies in a norm dense subset of $Y\,$.

Let $M$ and $N$ be von Neumann algebras.  The binormal tensor product
$N\otimes_{\t{bin}}M$ is the norm completion of $N\odot M$ with respect
to the norm
$$
\Vert x\Vert=\sup_\pi\Vert\pi(x)\Vert\quad,\quad x\in N\odot M ,
$$
where the supremum ranges over all representations of $N\odot M$
with normal restrictions to $N\otimes 1_M$ and $1_N\otimes M$,
called binormal representations \cite{EL}.

Let $F\subset(N\otimes_{\t{bin}}M)^*$ be the Banach space of normal
linear functionals associated with all such binormal representations.

\begin{lemma}\label{binmin} 
The kernel $J$ of the natural homomorphism
$$\pi:N\otimes_{\t{bin}}M\to N\otimes_{\min} M$$
is $\sigma(N\otimes_{\t{bin}}M,F)-\t{closed}$.
\end{lemma}

\proof When $M$ and $N$ act on $L^2(M)$ and $L^2(N)$ (respectively),
the $C^*$-subalgebra of $B(L^2(N)\overline\otimes L^2(M))$ generated by
$N\otimes 1_M$ and $1_N\otimes M$ carries the minimal tensor product norm
$\Vert\cdot\Vert_{\min}\,$. Thus we may identify it with
$N\otimes_{\min} M\,$.

With this identification, $\pi$ is a binormal representation of
$N\otimes_{\t{bin}}M$ on $L^2(N)\overline\otimes L^2(M)$,
whose kernel is $J$. Note now that $\pi$ is continuous from the 
$\sigma(N\otimes_{\t{bin}}M,F)$-topology to the $\sigma$-weak topology of
$B(L^2(N)\otimes L^2(M))$, therefore $J$ is
$\sigma(N\otimes_{\t{bin}}M,F)$-closed.
\endproof

\begin{proposition} Let $\alpha:G\to\t{Aut}(M)$ be an action of a
compact group $G$ on a von Neumann algebra $M$. Then $M$ is injective iff
the fixed point subalgebra $M^\alpha$ is injective.
\end{proposition}

\proof If $M$ is injective, so is $M^\a$
because there is a conditional expectation from $M$ to $M^\alpha$.

Now assume  that $M^\alpha$ is injective. In analogy with the proof of
Proposition \ref{stably}, we shall prove that, given any von Neumann
algebra $N\,$, the ideal $J$ in Lemma \ref{binmin} is $\{0\}$, what
means that $M$ is semidiscrete, that is injective \cite{CE,C2}.

The action $\beta=\iota\otimes \alpha$ on $N\odot M$ preserves the norm
$\Vert\cdot\Vert_{\t{bin}}\,$, so each $\beta_g$ extends to a
$*$-automorphism of $N\otimes_{\t{bin}}M$, still denoted by $\beta_g\,$.
Furthermore, the map $G\ni g\mapsto\beta_g(x)$ is
$\sigma(N\odot M,F)$-continuous for every $x\in N\odot M\,$. Therefore
the action $\beta$ on $N\otimes_{\t{bin}}M$ is
$\sigma(N\otimes_{\t{bin}}M,F)$-continuous.

Let us consider the conditional expectation $\E_\alpha:M\to M^\alpha$
defined by
$$\E_\alpha(x)=\int\alpha_g(x)\d g\ ,\qquad x\in M\ .$$
The completely positive map $\i\otimes \E_\a:N\otimes_{\min} M\to
N\otimes_{\min} M$ has norm 1, so
$$
\Vert(\i\otimes \E_\a)(x)\Vert_{\min}\leq\Vert x\Vert_{\min}\leq\Vert
x\Vert_{\t{bin}}\ ,\qquad x\in N\odot M\ .$$
Now let $x\in N\odot M$. Then $(\i\otimes \E_\a)(x)\in N\odot M^\alpha$. 
As $M^\alpha$ is semidiscrete, 
$\Vert\cdot\Vert_{\min}$ and $\Vert\cdot\Vert_{\t{bin}}$ 
coincide on $N\odot M^\alpha$. Therefore
$$
\Vert(\i\otimes \E_\a)(x)\Vert_{\t{bin}}\leq\Vert x\Vert_{\t{bin}}\ ,
\qquad x\in N\odot M,$$
so $(\i\otimes \E_\a)|_{N\odot M}$ extends to a  linear 
map of norm 1 from $N\otimes_{\t{bin}} M$ to itself, still denoted by 
$\i\otimes\E_\a$, whose
range is contained in the closure of $N\odot M^\alpha$ in 
$N\otimes_{\text{bin}}M$.

Let $x\in N\otimes_{\t{bin}}M$ and choose $x_n\in N\odot M$, $n\geq 1$, with
$\Vert x-x_n\Vert_{\t{bin}}\to 0$. Then, for all $\varphi\in F$,
\begin{multline*}
\varphi((\i_N\otimes\E_\alpha)(x)) = \lim_n\varphi((\i_N\otimes\E_\alpha)(x_n))\\
=\lim_n\int\varphi(\beta_g(x_n))\d g = \int\varphi(\beta_g(x))\d g\, ,
\end{multline*}
so
$$
\E_\beta(x) = \int\beta_g(x)\d g
$$
exists in the $\sigma(N\otimes_{\t{bin}}M,F)$-weak sense and 
$\E_\beta(x)=(\i_N\otimes\E_\alpha)(x)$ for all $x\in N\odot M$. 
In particular, the range of $\E_\beta$
is contained in the closure of $N\odot M^\alpha$. Moreover $\E_\beta$ is 
faithful.

Now let $x\in J$ be a positive element. Since each $\beta_g$ leaves $J$ 
globally invariant,
$\beta_g(x)\in J$ for all $g\in G$.
By the Hahn-Banach theorem $\E_\beta(x)$ belongs to the 
$\sigma(N\otimes_{\t{bin}}M,F)$-closed linear span of 
$\{\beta_g(x);\,g\in G\}$, thus
$\E_\beta(x)\in J$ by Lemma \ref{binmin}.

But $\E_\beta(x)$ is also in the closure of $N\odot M^\alpha$ which,
by the semidiscreteness of $M^\alpha$, intersects $J$ only in $0$.
Consequently $\E_\beta(x)=0$ and the faithfulness of $\E_\beta$ yields
$x=0$. We conclude that $J$ is trivial and $M$ is semidiscrete, thus 
injective. 
\endproof
\noindent
{\bf Acknowledgements.} We thank S. Woronowicz for a clarifying answer
to a question and F. P. Boca for calling our attention to
reference \cite{B}. 

\end{document}